\documentstyle[12pt]{amsart}
\newtheorem{th}{Theorem}[section]
\newtheorem{crl}[th]{Corollary}
\newtheorem{prp}[th]{Proposition}
\newtheorem{lm}[th]{Lemma}

\newcommand{\der }{\partial}

\newcommand{\m}{\mathop{\circ}\limits_{-}}
\newcommand{\limline}{\limits_{-\!-\!-\!-\!-\!-}}
\newcommand{\limsim}{\limits_{\sim\!\sim\!\sim\!\sim\!\sim\!\sim}}
\newcommand{\limeq}{\limits_{=\!=\!=\!=\!=\!=}}
\newcommand{\limsimeq}{\limits_{\simeq\!\simeq\!\simeq\!\simeq\!\simeq\!\simeq}}

\newcommand{\limsi}{\limits_{\sim\; \sim\; \sim\; \sim}}

\begin{document}
\title{Novikov-Jordan algebras}
\author {Askar Dzhumadil'daev}
\address
{Institute of Mathematics, Almaty, Kazakhstan}
\address{S. Demirel University, Almaty, Kazakhstan}
 \email{askar@@math.kz}

\maketitle

\begin{abstract} Algebras with identity $(a\star  b)\star  (c\star  d)
-(a\star d)\star(c\star  b)$ $=(a,b,c)\star  d-(a,d,c)\star  b$
are studied. Novikov algebras under Jordan multiplication and Leibniz
dual algebras  satisfy this identity. If algebra with a such identity has
unit, then it is associative and commutative.
\end{abstract}

\section{Introduction}
Let $(A,\circ)$ be some algebra with multiplication $A\times A\rightarrow A,$
$(a,b)\mapsto a\circ b.$  One can consider an algebra $A^+=(A,\{\;,\,\})$
under Jordan product $\{a,b\}=a\circ b+b\circ a.$ The algebra $A^+$ satisfies
the commutativity identity
$$\{a,b\}=\{b,a\}.$$

Let $\mathcal V$ be some variety of algebras. Denote by ${\mathcal V}^+$
a variety generated by algebras $A^+,$ where $A\in \mathcal V.$
One can ask about a minimal identity of minimal degree for $\mathcal V^+$
that does not follow from the commutativity identity.
If such identity exists call it {\it Jordan identity}.

Let $Ass$ be a class of associative algebras.
Well known, that Jordan identity for $Ass^+$ is the following identity
of degree 4
$$\{\{a,b\},\{c,d\}\}+
\{\{a,c\},\{d,b\}\}+
\{\{a,d\},\{b,c\}\}$$
$$-\{\{b,c\},a\},d\}\}-
\{\{c,d\},a\},b\}\}-
\{\{d,b\},a\},c\}=0.$$
An algebra $(A,\circ)$ with identity
$$ a\circ(b\circ c)-(a\circ b)\circ c=a\circ(c\circ
b)-(a\circ c)\circ b $$
is called {\it right-symmetric}.
This identity can be written in terms of
associators $(a,b,c)$ $= a\circ (b\circ  c)$ $-(a\circ  b)\circ
c,$ as follows $$(a,b,c)=(a,c,b).$$ That is the reason why this
identity is called right-symmetric. Right-symmetric algebras was
introduced around 1960's in \cite{Vinberg}, \cite{Kozshul},
\cite{Gerst1}.

A right-symmetric algebra $A$ is called (right) {\it Novikov}
\cite{nov}, \cite{Osborn1}, if the following {\it
left-commutativity} identity is true $$ a\circ(b\circ c)=b\circ
(a\circ c). $$

In fact right-symmetric was considered about 100 years before in
\cite{cayley} and Novikov algebras appeared earlier in
\cite{GelfandDorfman}. Right-symmetric algebras are
Lie-admissible. Any right-symmetric algebra under commutator
$[a,b]=a\circ  b-b\circ a$ is a Lie algebra.

An algebra with multiplication $A\times A\rightarrow A,$ $
(a,b)\mapsto a\times b$ and the identity
\begin{equation}
\label{Leibnizdual} (a\times b)\times c=a\times (b\times c+c\times
b)\end{equation} is called (left) Leibniz dual
\cite{Lodaycupprod}. Notice that any Leibniz dual algebra is
right-commutative
\begin{equation}\label{rcommut}(a\times b)\times
c=(a\times c)\times b.\end{equation}

For these algebras one can define their right and left versions by
considering algebras under the opposite multiplication
$(a,b)\mapsto b\circ a.$ In our paper we mainly consider  right
symmetric and right Novikov algebras. If otherwise are not stated
the names Leibniz dual and Novikov mean left Leibniz dual and right
Novikov, correspondingly.

We consider algebras with the identity of degree 4
\begin{equation}
\label{main}(a\star  b)\star  (c\star  d)-(a\star  d)\star (c\star
b)= (a,b,c)\star  d-(a,d,c)\star  b.
\end{equation}
Here $(a,b)\mapsto a\star  b$ is the multiplication of algebra and
$(a,b,c)$ $=a\star  (b\star  c)$ $-(a\star  b)\star  c$ is the
associator for the multiplication $\star$. Call this identity
a (right) {\it Tortken} identity. Call any algebra with
identity (\ref{main}) (right) {\it Tortken} algebra. In kazak
``tort'' or ``tirt'' means four.

Examples of Tortken algebras.

1. Polynomial algebra $K[x]$ under multiplication $$(a,b)\mapsto \der(ab).$$

2. Polynomial algebra $K[x]$ under multiplication $$(a,b)\mapsto
a\int_0^x b\,dx.$$

3. Divided power algebra $O_1(m)$ under multiplication
$$(a,b)\mapsto \der(\der^{p-1}(a)\der^{p-1}(b)),$$ if $char\,K=p>0.$

4. Divided power algebra $O_1(m)$ under multiplication $$(a,b)\mapsto
\der^{2^k+1}(\der^{2^k-1}(a)\der^{2^k-1}(b)),$$ if $p=2$ and $k>0.$

Let $\mathcal N$ be a class of Novikov algebras. Call $A^+$ {\it
Novikov-Jordan} , if $A$ is Novikov. We show that Jordan identity
for $\mathcal N^+$ is Tortken identity. We prove that some Osborn
algebras under Jordan product are simple. We establish that
Leibniz dual algebras gives us another examples of algebras with
identity (\ref{main}). These algebras are not commutative, but
they are right-commutative. So, in general the class of Tortken
algebras does not coincide with the class of Novikov-Jordan
algebras.

One can consider left Tortken identity
$$(a\star b)\star (c\star d)-(c\star b)\star (a\star d)= -a\star
(b,c,d)+c\star (b,a,d).$$ If $(A,\star)$ is right Tortken, then
$(A, \bar \star),$  where $a\bar\star b=b\star a$ is left Tortken. For
commutative Tortken algebras, like Novikov-Jordan algebras, left
and right Tortken identities are coincide. In general they are
different identities. 

We prove  that algebras with identity (\ref{main}) and with a unit
is associative and commutative. In \cite{os1}, \cite{os2},
\cite{os3}, \cite{os4}, \cite{os5} algebras with identity of
degree 3 and 4 are classified. Their classification was done under
the condition that algebras have unit. Therefore, the class of
algebras with identity (\ref{main}) is not in  their list.

Let $[a,b,c]=(a,b,c)-(a,c,b)$ be a deviation of associators.  In
\cite{sokolov} is proved that many known classes of algebras (Lie,
Jordan, right-symmetric algebras, LT-algebras) are satisfy the
identity
$$[a\ast b,c,d]=a\ast[b,c,d]+[a,c,d]\ast b.$$
Novikov-Jordan algebras $Os^+(\alpha,\beta)$ do not satisfy this
identity. Therefore, algebras with identity (\ref{main}) are not
in this class.

Identities of right-symmetric algebras was also studied in \cite{Dzh00}.

\section{Tortken identity and some consequences of degree 4}

Let $A$ be an algebra over the field $K$ of characteristic $p\ge
0$ and $f=f(t_1,\ldots,t_k)$ be some non-associative polynomial in
variables $t_1,\ldots, t_k.$ We say that $f=0$ or
$f(a_1,\ldots,a_k)=0$ is identity on $A,$ if
$f(a_1,\ldots,a_k)=0$ for any substitution $t_i:=a_i\in A.$

Let $A$ be an algebra $A$ with multiplication $A\times
A\rightarrow A, (a,b)\mapsto a\star b.$ Recall that its element
$e\in A$ is called a {\it left unit,} if $e\star a=a,$ for any
$a\in A.$ Similarly, $e$ is called a {\it right unit,} if $a\star
e=a,$ for any $a\in A.$  Left and right unit is called {\it unit.}

\begin{th}\label{unit}
Let $A$ be an algebra with the identity $(\ref{main})$.

{\rm i)} If $A$ has right unit, then $A$ satisfies the  following
identity
\begin{equation}
\label{rightunit} a\star (b\star c+c\star b)=2 (a\star b)\star c.
\end{equation}

{\rm ii)} If $A$ has left unit, then $A$ is associative and commutative.
\end{th}

{\bf Proof.} i) Let $e$ be the right unit. Then
\begin{equation}\label{komekshi2}
(a,b,e)=(a\star b)\star e-a\star (b\star e)=0,
\end{equation}
for any $a,b,c\in A.$ Take $c:=e$ in (\ref{main}). By
(\ref{komekshi2}) we have
\begin{equation}
\label{komekshi1} (a\star b)\star (e\star d)-(a\star d)\star
(e\star b)=0,
\end{equation}
for any $a,b,d\in A.$ In (\ref{komekshi1}) take $d:=e.$ We obtain
that
$$a\star b-a\star (e\star b)=$$
$$(a\star b)\star(e\star e)-(a\star e)\star (e\star b)=0,$$
for any $a, b\in A.$ In other words,
\begin{equation}
\label{komekshi} (a\star e)\star b-a\star (e\star b)=(a,e,b)=0.
\end{equation}

Now take $d:=e$ in (\ref{main}). By (\ref{komekshi}) we have
$$
(a\star  b)\star  (c\star  e)
-(a\star  e)\star  (c\star  b)=(a,b,c)\star  e.
$$
So,
$$(a\star  b)\star  c-a\star  (c\star  b)=(a,b,c).$$
In other words, (\ref{rightunit}) is true.

ii) Let $e$ be the left unit. Then $$(e,a,b)=e\star  (a\star
b)-(e\star  a)\star  b=0,$$ for any $a,b\in A.$ Therefore, for
$a:=e$ from (\ref{main}) it follows that
\begin{equation}\label{7nov}
b\star  (c\star
d)-d\star (c\star  b)=0.
\end{equation}
In (\ref{7nov}) take $c:=e.$ We
have $$b\star  d-d\star  b=0,$$ for any $b,d\in A.$ So, the
algebra $A$ is commutative. Thus, according (\ref{7nov}),
$$b\star (c\star d)=d\star(c\star b)=(c\star b)\star d=(b\star c)\star d.$$

\begin{crl} \label{unit1} Any algebra with unit and identity $(\ref{main})$ is
associative and commutative.
\end{crl}

{\bf Remark.} If $p\ne 2$ from the identity (\ref{rightunit})
follow the right-alternative identity $(a,b,b)=0$ and
right-commutativity identity $(a\star  b)\star  c=(a\star
c)\star  b$ and vice versa, from these two identities follow
(\ref{rightunit}). Notice that, the right-alternativity identity
itself is not enough to obtain (\ref{rightunit}).

For an algebra $A$ with multiplication $A\times A\rightarrow A,
(a,b)\mapsto a\star b,$ let $r_a: A\rightarrow A, b\mapsto a\star
b$ be the operator of multiplication to $a.$

\begin{prp}\label{rsym3} Let $A$ be a commutative algebra with identity
$(\ref{main})$. Then $$\sum_{\sigma\in Sym_3}sign\,\sigma\,
r_{a_{\sigma(1)}}r_{a_{\sigma(2)}}r_{a_{\sigma(3)}}=0,$$ for any
$a_1,a_2,a_3\in A.$
\end{prp}

{\bf Proof.}
$$((x\star a)\star b)\star c+ ((x\star b)\star c)\star a+
((x\star c)\star a)\star b$$
$$-((x\star a)\star c)\star b-((x\star b)\star a)\star c- ((x\star c)\star b)\star a=$$

$$((x\star a)\star b)\star c-((x\star b)\star a)\star c+
\mathop{((c\star b)\star a)\star x}\limline -\mathop{((c\star
a)\star b)\star x}\limsim+$$
$$((x\star b)\star c)\star a-((x\star c)\star b)\star a-
\mathop{((a\star b)\star c)\star x}\limeq+\mathop{((a\star c)\star
b)\star x}\limsim+$$
$$((x\star c)\star a)\star b)-((x\star a)\star c)\star b+
\mathop{((b\star a)\star c)\star x}\limeq -\mathop{((b\star
c)\star a)\star x}\limline=$$ (by identity (\ref{main}))
$$((b\star x)\star (a\star c))-((b\star c)\star (a\star x))+$$
$$((c\star x)\star (b\star a))-((c\star a)\star (b\star x))+$$
$$((a\star x)\star (c\star b))-((a\star b)\star (c\star x))$$ (by the
commutativity condition) $$=0.$$

\begin{crl} \label{on} Let $A$ be a commutative algebra with identity
$(\ref{main})$. Then for any $a,b,c,x\in A,$
$$(a,b\star x,c)+(b,c\star x,a)+(c,a\star x,b)=0. $$
\end{crl}

{\bf Proof.} By proposition~\ref{rsym3} $$(c,a\star x,b)+(a,b\star
x,c)+(b,c\star x,a)=$$ $$((x\star a)\star b)\star c -((x\star
a)\star c)\star b)$$ $$+((x\star b)\star c)\star a)-((x\star
b)\star a)\star c)$$ $$+((x\star c)\star a)\star b) -((x\star
c)\star b)\star a)$$ $$=0.$$

\begin{crl}
\label{onon} Let $A$ be a commutative algebra with identity
$(\ref{main})$. Then for any $a,b,c,x\in A,$
$$ (a,x,b)\star c+(b,x,c)\star a+(c,x,a)\star b=0.$$
\end{crl}

{\bf Proof.} By proposition~\ref{rsym3},
$$(a,x,b)\star c+(b,x,c)\star a+(c,x,a)\star b=$$
$$((x\star b)\star a)\star c)-((x\star a)\star b)\star c)$$
$$+((x\star c)\star b)\star a)-((x\star b)\star c)\star a)$$
$$+((x\star a)\star c)\star b)-((x\star c)\star a)\star b)=$$
$$=0.$$

\section{Novikov algebras under Jordan product\label{nov-jor}}

\begin{prp}\label{birbir} Let $A$ be a Novikov algebra. Then for any
$a_1,a_2,a_3\in A,$ $$\sum_{\sigma\in Sym_3}sign\,\sigma\,
r_{a_{\sigma(1)}}r_{a_{\sigma(2)}} r_{a_{\sigma(3)}}=0.$$ Here
$r_a$ be the right-multiplication operator: $(b)\, r_a=b\circ a.$
\end{prp}

{\bf Proof.} For any $a,b,c,d\in A,$ we have
$$ ((a\circ b)\circ
c)\circ d+((a\circ c)\circ d)\circ b +((a\circ d)\circ b)\circ c$$
$$-((a\circ b)\circ d)\circ c-((a\circ c)\circ b)\circ d -((a\circ
d)\circ c)\circ b=$$ (by right-symmetric rule) $$(a\circ b)\circ
(c\circ d-d\circ c)+(a\circ c)\circ (d\circ b-b\circ d) +(a\circ
d)\circ (b\circ c-c\circ b)=$$ (by left-commutativity rule)
$$c\circ((a\circ b)\circ d)- d\circ((a\circ b)\circ c)$$ $$+d\circ
((a\circ c)\circ b) -b\circ ((a\circ c)\circ d)$$ $$+b\circ
((a\circ d)\circ c)- c\circ ((a\circ d)\circ b)=$$ (by
right-symmetric rule) $$b\circ (a\circ (d\circ c-c\circ d))$$
$$+c\circ(a\circ (b\circ d-d\circ b))$$ $$+d\circ(a\circ (c\circ
b-b\circ c))=$$ (by left-commutativity rule) $$b\circ (d\circ
(a\circ c))-b\circ(c\circ(a\circ d))$$ $$+c\circ(b\circ (a\circ
d))-c\circ(d\circ(a\circ b))$$ $$+d\circ(c\circ (a\circ
b))-d\circ(b\circ(a\circ c))$$ (by left-commutativity rule)
$$=0.$$

\begin{crl} Let $a\m b= a\circ f(b)$ be
a new multiplication in Novikov algebra $(A,\circ),$ where $f\in
End\,A.$ Then the algebra $(A,\m)$ satisfies the following
identity of degree $4$ $$((a\m b)\m c)\m d+((a\m c)\m d)\m b+((a\m
d)\m b)\m c$$ $$-((a\m c)\m b)\m d-((a\m d)\m c)\m b-((a\m b)\m
d)\m c=0,$$ for any $a,b,c\in A.$
\end{crl}

{\bf Proof.} Notice that $a\m b=ar_{f(b)}.$ Therefore our identity
is equivalent to the relation $$\sum_{\sigma\in
Sym_3}sign\,\sigma\, r_{f(a_{\sigma(1)})}r_{f(a_{\sigma(2)})}
r_{f(a_{\sigma(3)})}=0.$$

{\bf Example 1.} Any associative algebra is right-symmetric.

{\bf Example 2.} Let $U$ be the algebra of Laurent polynomials
$K[x^{\pm 1}]$ (if $p=0$) or divided power algebra $O_1(m)$ (if
$p>0$) and $W_1$ or $W_1(m)$ be the Lie algebra of special
derivations $u\der,$ where $\der(x^i)=ix^{i-1}$ (if $p=0$) or
$\der(x^{(i)})=x^{(i-1)}.$ Define the multiplication on $W_1$ or
$W_1(m)$ by $$u\der\circ v\der=v\der(u)\der.$$ We obtain a
right-symmetric algebra. Call it as a right-symmetric Witt
algebra. This construction can be easily generalized for the case
of any associative commutative algebra $U$ with derivation $\der.$

{\bf Example 3.} An algebra $A$ with a basis $\{e_1,\ldots,e_n\}$
and multiplication $e_i\circ e_j=e_j$ is left-commutative and
associative. In particular, $A$ is Novikov. Its Jordan algebra
arises in genetics and it is called as a gametic algebra
\cite{holgate}. Jordan algebra for $A$ satisfies
(associative)-Jordan identity
$\{\{\{x,x\},y\},x\}=\{\{x,x\},\{y,x\}\}$ and identity
(\ref{main}).

Any associative algebra is Lie-admissible: under commutator
$[a,b]=a\circ b-b\circ a$ it can be endowed by a structure of Lie
algebra. It is also well known, that an associative algebra with
multiplication $(a,b)\mapsto a\circ b$ can by endowed by a
structure of Jordan algebra under the anti-commutator
$\{a,b\}=a\circ b+b\circ a.$ Right-symmetric algebras are also
Lie-admissible. Below we describe analog of Jordan algebras for
Novikov algebras.

A commutative Tortken algebra is called {\it special Tortken,}
if it is isomorphic to some subalgebra of Tortken algebra of the form
$A^+,$ where $A$ is Novikov.

\begin{th}\label{one}
Let $A$ be a Novikov algebra. Then the algebra $A^+$ satisfies the
Tortken identity; moreover, this identity is a Jordan identity for
$A^+.$
\end{th}

Proof of theorem \ref{one} needs several lemmas. Below we assume
that $A$ is Novikov algebra and that $a,b,c,d\in A.$

\begin{lm} \label{four}
$$(a,b,c)^+=\{\{b,c\},a\}-\{\{b,a\},c\}= \{b,[c,a]\}.$$
\end{lm}

{\bf Proof.} It is easy to see that
 $$(a,b,c)^+=\{\{b,c\},a\}-\{\{b,a\},c\}.$$
Therefore, $$(a,b,c)^+=$$
 $$(b\circ c)\circ a+(c\circ b)\circ
a+a\circ(b\circ c)+ \mathop{a\circ(c\circ b)}\limline$$
$$-(b\circ a)\circ c-(a\circ b)\circ c-c\circ(b\circ a)-
\mathop{c\circ(a\circ b)}\limline=$$
(left-commutativity identity)
$$=(b\circ c)\circ a-(b\circ a)\circ c -(a\circ c)\circ
b+(c\circ a)\circ b+\mathop{(a\circ c)\circ b}\limeq$$
$$-\mathop{(c\circ a)\circ b}\limsi+ \mathop{(c\circ b)\circ
a}\limsi -\mathop{(a\circ b)\circ c}\limeq +a\circ (b\circ
c)-c\circ(b\circ a).$$
By left-commutativity,
$$a\circ (b\circ c)-c\circ(b\circ a)=$$
$$b\circ (a\circ c)-b\circ (c\circ a)=$$
$$b\circ [a,c].$$
Therefore,
$$(a,b,c)^+=$$
$$b\circ [c,a]+[c,a]\circ b+a\circ [c,b]+c\circ [b,a]+b\circ[a,c]=$$
$$\{b,[c,a]\}.$$

\begin{lm}\label{five}
$$\{(a,b,c)^+,d\}-\{(a,d,c)^+,b\}=
\{[c,a],[b,d]\}.$$
\end{lm}

{\bf Proof.} By lemma \ref{four},
$$\{(a,b,c)^+,d\}=\{\{b,[c,a]\}, d\},$$
$$\{(a,d,c)^+,b\}= \{\{d,[c,a]\},b\}.$$
Thus,
$$\{(a,b,c)^+,d\}-\{(a,d,c)^+,b\}= $$
$$\{\{b,[c,a]\}, d\}- \{\{d,[c,a]\},b\}=$$
$$(d,[c,a],b)^+=$$
$$\{[c,a],[b,d]\}.$$

\begin{lm} \label{six}
$$\{\{a,b\},\{c,d\}\}-\{\{a,d\},\{c,b\}\}=$$
$$\{d\circ c,a\circ b\}+\{c\circ d,b\circ a\}-
\{b\circ c,a\circ d\}-\{c\circ b,d\circ a\}.$$
\end{lm}

{\bf Proof.} We have
$$\{\{a,b\},\{c,d\}\}-\{\{a,d\},\{c,b\}\}=$$
$$(a\circ b)\circ (c\circ d)+(b\circ a)\circ (c\circ d)+
(a\circ b)\circ (d\circ c)+(b\circ a)\circ (d\circ c)$$
$$+(c\circ d)\circ (a\circ b)+(c\circ d)\circ (b\circ a)+
(d\circ c)\circ (a\circ b)+(d\circ c)\circ (b\circ a)$$
$$-(a\circ d)\circ (c\circ b)-(a\circ d)\circ (b\circ c)
-(d\circ a)\circ (c\circ b)-(d\circ a)\circ (b\circ c)$$
$$-(c\circ b)\circ (a\circ d)-(b\circ c)\circ (a\circ d)
-(c\circ b)\circ (d\circ a)-(b\circ c)\circ (d\circ a)=$$
(left-commutativity rule)
$$=c\circ((a\circ b)\circ d)+c\circ((b\circ a)\circ d)+
d\circ((a\circ b)\circ c)+d\circ((b\circ a)\circ c)$$
$$+a\circ((c\circ d)\circ b)+b\circ((c\circ d)\circ a)+
a\circ((d\circ c)\circ b)+b\circ((d\circ c)\circ a)$$
$$-c\circ((a\circ d)\circ b)-b\circ((a\circ d)\circ c)
-c\circ((d\circ a)\circ b)-b\circ((d\circ a)\circ c)$$
$$-a\circ((c\circ b)\circ d)-a\circ((b\circ c)\circ d)
-d\circ((c\circ b)\circ a)-d\circ((b\circ c)\circ a)=$$
$$+a\circ((c\circ d)\circ b+(d\circ c)\circ b
-(c\circ b)\circ d-(b\circ c)\circ d)$$
$$+b\circ((c\circ d)\circ a+(d\circ c)\circ a
-(a\circ d)\circ c-(d\circ a)\circ c)$$
$$+c\circ((a\circ b)\circ d+(b\circ a)\circ d
-(a\circ d)\circ b-(d\circ a)\circ b)$$
$$+d\circ((a\circ b)\circ c)+(b\circ a)\circ c
-(c\circ b)\circ a-(b\circ c)\circ a)=$$
(right-symmetric and
left-commutativity rules )
$$=\mathop{a\circ (c\circ [d,b])}\limline+
(d\circ c)\circ (a\circ b)-(b\circ c)\circ (a\circ d)$$
$$+\mathop{b\circ (d\circ [c,a])}\limeq+
(c\circ d)\circ (b\circ a)-(a\circ d)\circ (b\circ c)$$
$$+\mathop{c\circ (a\circ [b,d])}\limline+
(b\circ a)\circ (c\circ d)-(d\circ a)\circ (c\circ b)$$
$$+\mathop{d\circ (b\circ [a,c])}\limeq+
(a\circ b)\circ (d\circ c)-(c\circ b)\circ (d\circ a)=$$
(left-commutativity rule)
$$(d\circ c)\circ (a\circ b)-(b\circ c)\circ (a\circ d)+
(c\circ d)\circ (b\circ a)-(a\circ d)\circ (b\circ c)$$
$$+(b\circ a)\circ (c\circ d)-(d\circ a)\circ (c\circ b)+
(a\circ b)\circ (d\circ c)-(c\circ b)\circ (d\circ a)=$$
$$\{d\circ c,a\circ b\}+\{c\circ d,b\circ a\}-
\{b\circ c,a\circ d\}-\{d\circ a,c\circ b\}.$$

\begin{lm}\label{seven}
$$\{[c,a],[b,d]\}=$$
$$\{d\circ c,a\circ b\}+\{c\circ d,b\circ a\}-
\{b\circ c,a\circ d\}-\{c\circ b,d\circ a\}.$$
\end{lm}

{\bf Proof.} We have
$$\{[c,a],[b,d]\}-\{d\circ c,a\circ b\}-\{c\circ d,b\circ a\}+
\{b\circ c,a\circ d\}+\{c\circ b,d\circ a\}=$$
$$(c\circ a)\circ (b\circ d)-(c\circ a)\circ (d\circ b)-
(a\circ c)\circ (b\circ d)+(a\circ c)\circ (d\circ b)$$
$$+(b\circ d)\circ (c\circ a)-(d\circ b)\circ (c\circ a)-
(b\circ d)\circ (a\circ c)+(d\circ b)\circ (a\circ c)$$
$$-(d\circ c)\circ (a\circ b)-(a\circ b)\circ (d\circ c)$$
$$-(c\circ d)\circ(b\circ a)-(b\circ a)\circ (c\circ d)$$
$$+(b\circ c)\circ(a\circ d)+(a\circ d)\circ (b\circ c)$$
$$+(c\circ b)\circ (d\circ a)+(d\circ a)\circ (c\circ b)=$$
(left-commutativity rule)
$$=b\circ((c\circ a)\circ d)-d\circ((c\circ a)\circ b)-
b\circ((a\circ c)\circ d)+d\circ((a\circ c)\circ b)$$
$$+c\circ((b\circ d)\circ a)-c\circ((d\circ b)\circ a)-
a\circ((b\circ d)\circ c)+a\circ((d\circ b)\circ c)$$
$$-a\circ((d\circ c)\circ b)-d\circ((a\circ b)\circ c)$$
$$-b\circ((c\circ d)\circ a)-c\circ((b\circ a)\circ d)$$
$$+a\circ((b\circ c)\circ d)+b\circ((a\circ d)\circ c)$$
$$+d\circ((c\circ b)\circ a)+c\circ((d\circ a)\circ b)=$$
$$a\circ(-(b\circ d)\circ c+(d\circ b)\circ c
-(d\circ c)\circ b+(b\circ c)\circ d)$$
$$+b\circ((c\circ a)\circ d-(a\circ c)\circ d
-(c\circ d)\circ a+(a\circ d)\circ c)$$
$$+c\circ((b\circ d)\circ a-(d\circ b)\circ a-
(b\circ a)\circ d+(d\circ a)\circ b)$$
$$+d\circ(-(c\circ a)\circ b)+(a\circ c)\circ b
-(a\circ b)\circ c+(c\circ b)\circ a)=$$
$$a\circ(b\circ [c,d]+d\circ [b,c])$$
$$+b\circ(c\circ [a,d]-a\circ [c,d])$$
$$+c\circ (b\circ [d,a]-d\circ [b,a])$$
$$+d\circ (-c\circ [a,b]+a\circ [c,b])=$$
(left-commutativity rule)
$$a\circ(b\circ [c,d])+a\circ(d\circ [b,c])$$
$$+b\circ(c\circ [a,d])-a\circ(b\circ [c,d])$$
$$+b\circ (c\circ [d,a])-c\circ(d\circ [b,a])$$
$$+c\circ (-d\circ [a,b])+a\circ(d\circ [c,b])=$$
$$0.$$
Lemma \ref{seven} is proved.

Denote by $O_1(\infty)$ infinite-dimensional algebra $\{x^{(i)}: 0\le i\}$
with multiplication
$$x^{(i)} x^{(j)}={i+j\choose i} x^{(i+j)}.$$
If $p>0,$ then it has $p^m$-dimensional subalgebra
$$O_1(m)=\{x^{(i)}: 0\le i<p^m\},$$
for any nonnegative integer $m.$ The algebra $O_1(m)$ is called
{it divided power algebra.} If $p=0,$ then $O_1(\infty)$ is
isomorphic to polynomial algebra $K[x].$ The isomorphism can be
given by $x^{(i)}\mapsto x^i/i!.$ Take now a new basis of the
vector space $O_1(\infty)$ that consists of elements
$e_i=x^{(i+1)},$ $-1\le i.$ Endow the vector space  $O_1(m),$
$p>0$ by multiplication
$$e_i\star e_j={i+j+2\choose i+1}e_{i+j}.$$
This algebra is commutative and satisfies identity (\ref{main}).
In notation of section \ref{osborn} this algebra is isomorphic to
$Os^+(0,0,m).$

\begin{lm} \label{16n} $(p\ne 2).$
Let $Os^+(0,0)=\{e_i: -1\le i , i\in {\bf Z}\}$ be Novikov-Jordan
algebra with multiplication $e_i\star e_j=(i+j+2)e_{i+j}$ (we will
consider this algebra more detailed in section $(\ref{osborn})$).
Any multilinear polynomial identity of degree $3$ of the algebra
$(Os^+(0,0),\star)$ follows from commutativity identity.
\end{lm}

{\bf Proof.} According commutativity rule we can assume that $f$ has the form
$$f(t_1,t_2,t_3)=
\lambda_1(t_1t_2)t_3+\lambda_2(t_2t_3)t_1+\lambda_3(t_3t_1)t_2.$$
Then
$$f(e_i,e_j,e_s)=
((i+j+s+2)(i+j+2)\lambda_1+(j+s+2)\lambda_2+(s+i+2)\lambda_3)e_{i+j+1}.$$
Therefore, the condition $f(e_i,e_j,e_s)=0$ gives us that
$$(i+j+2)\lambda_1+(j+s+2)\lambda_2+(s+i+2)\lambda_3=0,$$
for any $i,j,s\in {\bf Z}, i,j,s\ge -1,$ such that $i+j+s+2\ne 0.$
So, for $(i,j,s)=(1,2,3),(2,3,1),(3,1,2),$ we obtain the system of
linear equations, with determinant $54.$ Thus, this system is
non-degenerate, if $p\ne 2,3.$ Let $p=3.$ Then for
$(i,j,s)=(1,1,0), (1,0,1), (0,1,1)$ we obtain non-degenerate
system of $3\times 3$ equations. Hence,
$\lambda_1=\lambda_2=\lambda_3=0,$ if $p\ne 2.$

{\bf Remark.} If $p=2,$ lemma \ref{16n} is not true. In this case Jordan
product and Lie products are coincide and any special Tortken
algebra satisfies Jacobi identity.

\begin{lm} \label{16november} $(p\ne 2).$
For the algebra $(O_1(m),\star)$ the  space of
multilinear polynomial identities of
degree $4$ is generated  by the following polynomials
$$Tortken(t_1,t_2,t_3,t_4):=$$
$$(t_1t_2)(t_3t_4)-(t_1 t_4)(t_3t_2)-(t_1(t_2t_3))t_4+((t_1t_2)t_3)t_4+
(t_1(t_4t_3))t_2-((t_1t_4)t_3)t_2,$$

$$Tortken'(t_1,t_2,t_3,t_4):=$$
$$(t_1t_3)(t_2t_4)+(t_1t_4)(t_2t_3)+((t_1t_3)t_4)t_2$$
$$+((t_1t_4)t_2)t_3+((t_2t_3)t_1)t_4+((t_2t_4)t_3)t_1, \mbox{\;if\;\;} (p,m)=(3,1),$$

$$Com(t_1,t_2):=t_1t_2-t_2t_1.$$
\end{lm}

{\bf Proof.}
According commutativity condition we can assume that
$f$ has the following form
$$f(t_1,t_2,t_3,t_4)=$$
$$\mu_1 (t_1t_2)(t_3t_4)+\mu_2(t_1t_3)(t_2t_4)+\mu_3(t_1t_4)(t_2t_3)$$
$$\mu_4((t_1t_2)t_3)t_4+\mu_5((t_1t_2)t_4)t_3+\mu_6((t_1t_3)t_2)t_4+
\mu_7((t_1t_3)t_4)t_2$$
$$\mu_8((t_1t_4)t_2)t_3+\mu_9((t_1t_4)t_3)t_2+\mu_{10}((t_2t_3)t_1)t_4+
\mu_{11}((t_2t_3)t_4)t_1$$
$$+\mu_{12}((t_2t_4)t_1)t_3+\mu_{13}((t_2t_4)t_3)t_1+\mu_{14}((t_3t_4)t_1)t_2+
\mu_{15}((t_3t_4)t_2)t_1.$$

\medskip
Suppose that $p>3$ or $m>1,$ if $p=3.$
Let $1=x^{(0)},  x=x^{(1)}.$
Make the following $10$ substitutions of $(t_1,t_2,t_3,t_4):$
$$(x,x,x,1), (x,x,1,x),(x,1,x,x),(1,x,x,x),
(1,1,x,x^{(2)}),(1,1,x^{(2)},x),$$
$$
(1,x^{(2)},x,1), (x,1,1,x^{(2)}), (x^{(2)},1,1,x),
(1,1,1,x^{(3)}).$$

\noindent We obtain the system of linear equations
$M{\mu^t}=0,$ where ${\mu}=(\mu_1,\ldots,\mu_{15})$ is a
row with $15$ components, $\mu^t$ is corresponding column and
$M$ is $10\times 15$-matrix
\bigskip

$$\left(\begin{array}{ccccccccccccccc}
2&2&2&4&2&4&2&1&1&4&2&1&1&1&1\\
2&2&2&2&4&1&1&4&2&1&1&4&2&1&1\\
2&2&2&1&1&2&4&2&4&1&1&1&1&4&2\\
2&2&2&1&1&1&1&1&1&2&4&2&4&2&4\\
0&1&1&0&0&0&1&1&2&0&1&1&2&3&3\\
0&1&1&0&0&1&2&0&1&1&2&0&1&3&3\\
1&1&0&2&1&1&0&0&0&3&3&1&2&0&1\\
1&1&0&0&1&0&1&3&3&0&0&2&1&2&1\\
1&1&0&1&2&1&2&3&3&0&0&1&0&1&0\\
0&0&0&0&0&0&0&1&1&0&0&1&1&1&1\\
\end{array}\right)
$$

\bigskip

\noindent It has rank $10$ and has the following fundamental
system of solutions
\bigskip

$$\begin{array}{ccccccccccccccc}
\{(1,&0,&-1,&0,&1,&0,&0,&-1,&0,&0,&-1,&0,&0,&0,&1),\\
(1,&0,&-1,&0,&1,&1,&-1,&-1,&0,&-1,&0,&0,&0,&1,&0),\\
(0,&1,&-1,&-1,&1,&1,&0,&-1,&0,&0,&-1,&0,&1,&0,&0),\\
(0,&1,&-1,&0,&0,&1,&0,&-1,&0,&-1,&0,&1,&0,&0,&0),\\
(0,&0,&0,&-1,&1,&1,&-1,&-1,&1,&0,&0,&0,&0,&0,&0)\}\\
\end{array}
$$

\bigskip

\noindent So, we can take $\mu_9,\mu_{12},\mu_{13},\mu_{14},\mu_{15}$
as a free parameters. Then

$\mu_1=\mu_{14} + \mu_{15},$

$\mu_2= \mu_{12} + \mu_{13},$

$\mu_3= -\mu_{12} -\mu_{13} - \mu_{14} - \mu_{15},$

$\mu_4= -\mu_{13} - \mu_9,$

$\mu_5 = \mu_{13} + \mu_{14} + \mu_{15} + \mu_9,$

$\mu_6 = \mu_{12} + \mu_{13} + \mu_{14} + \mu_9,$

$\mu_7 = -\mu_{14} - \mu_9,$

$\mu_8 = -\mu_{12} - \mu_{13} - \mu_{14} - \mu_{15} - \mu_9,$

$\mu_{10} = -\mu_{12} - \mu_{14},$

$\mu_{11} = -\mu_{13} - \mu_{15}.$

\noindent Thus,
$$f=\mu_9f_1+\mu_{12}f_2+\mu_{13}f_3+\mu_{14}f_4+\mu_{15}f_5,$$
where
$$f_1= -\sum_{\sigma\in Sym_3}sign\,\sigma\,((t_{1} t_{\sigma(2)})t_{\sigma(2)})t_{\sigma(3)},$$
$$f_2=Tortken(t_1,t_3,t_2,t_4),$$

$$f_3=(t_1t_3)(t_2t_4)-(t_1t_4)(t_2t_3)-((t_1t_2)t_3)t_4$$
$+((t_1t_2)t_4)t_3+((t_1t_3)t_2)t_4-((t_1t_4)t_2)t_3
-((t_2t_3)t_4)t_1+((t_2t_4)t_3)t_1,$

$$f_4=(t_1t_2)(t_3t_4)-(t_1t_4)(t_2t_3)+((t_1t_2)t_4)t_3 $$
$+ ((t_1t_3)t_2)t_4
- ((t_1t_3)t_4)t_2
- ((t_1t_4)t_2)t_3
- ((t_2t_3)t_1)t_4
+ ((t_3t_4)t_1)t_2,$

$$f_5=(t_1t_2)(t_3t_4)
- (t_1t_4)(t_2t_3) $$ $+ ((t_1t_2)t_4)t_3) - ((t_1t_4)t_2)t_3 -
((t_2t_3)t_4)t_1 + ((t_3t_4)t_2)t_1.$

\medskip

\noindent We see that, $f_2$  is the identity (\ref{main}) for
$t_1,t_3,t_2,t_4.$ By proposition~\ref{rsym3}, $f_1=0$ is identity
for any commutative Tortken algebra. It is easy to see that
$f_3=0, f_4=0, f_5=0$ are follow from corollary~\ref{onon}  and
(\ref{main}).

We omit details of tedious calculations for  $(p,m)=(3,1).$
They are similar to above given case. We just mention that
$$(a\star c)\star (b\star d)+(a\star d)\star (b\star c)
+((a\star c)\star d)\star b$$
$$+((a\star d)\star b)\star c+((b\star c)\star a)\star d+
((b\star d)\star c)\star a=2\der^3(a\cdot b\cdot c\cdot d),
$$
if $p=3$ and $a\star b=\der (a\cdot b).$ Thus
$Tortken'$ is an identity for $(p,m)=(3,1).$

{\bf Proof of theorem~\ref{one}.}
Tortken identity for Jordan product follows from lemma~\ref{five}, \ref{six}, \ref{seven}. Second part of theorem \ref{one} follows from
lemma~\ref{16n} and \ref{16november}.

\section{Simple Novikov-Jordan algebras\label{osborn}}

Let $O_1(m)=\{x^{(i)}: 0\le i\le p^m-1\}$ be divided power algebra with
multiplication
$$x^{(i)} x^{(j)}={i+j\choose j}x^{(i+j)},$$
if $p>0.$
If $p=0$ we can consider an algebra $O_1(\infty)=\{x^{(i)}: i\in {\bf Z}, i\ge 0\}$ with the same multiplication.

Let $U=K[[x^{\pm}]]$ be Laurent formal power series algebra, if $p=0,$
and $U=O_1(m),$ if $p>0.$  So, any element of $U$ for $p=0$ has a form
$\sum_{j\ge i} \lambda_jx^j,$  where a number of non-zero terms in
positive part $\lambda_j\in K, j\ge 0,$ may be infinite, but the number of
non-zero terms in negative part, $|\{\lambda_j: j<0\}|$  is finite.

In \cite{Osborn1}, \cite{Osborn2}, \cite{Osborn3}, it is proved
that the algebra defined on polynomials algebra $U=K[x]$ by the rule
$$(a,b)\mapsto \der(a)b$$
and on Laurent formal power series algebra
$U=K[[x^{\pm 1}]], p=0,$ by
$$(a,b)\mapsto \der(a)b+\alpha x^{-1} ab+\beta x^{-2}ab, \quad
\alpha,\beta\in K,$$ is Novikov and simple.
In the case of $p>0,
U=O_1(m)$ this multiplication is changed by obvious way
$$(a,b)\mapsto \der(a)b+\alpha x^{(p^m-1)} ab+\beta x^{(p^m-2)}ab,
\quad \alpha,\beta\in K.$$ Denote such algebras in honor of
M.~Osborn by $Os$ if $U=K[x],p=0,$ $Os(\alpha,\beta),$ if $U=K[[x^{\pm}]],
p=0,$ and $Os(\alpha,\beta,m),$ if $U=O_1(m), p>2.$

In this section we denote by  $\star$ the Jordan bracket $\{\;,\;\}.$

Let $p=0.$ Take  $U=K[[x^{\pm}]].$ Let
$\overline{Os}(\alpha,0)$ be a subspace of $Os(\alpha,0)$
generated by the set  $\{x^{i}: i\ne -2\alpha-1,i\in {\bf Z}\},$
if $\alpha\in {\frac{1}{2}\bf Z}.$ Notice that
$$x^i\star x^j=(i+j+2\alpha)x^{i+j-1}+2\beta
x^{i+j-2}.$$ Therefore, the coefficient at  $<x^{-2\alpha-1}>$
of any element of the form
$X=x^{i}\star \sum_{j}\lambda_jx^j=\sum_j\lambda_j(i+j+2\alpha)x^{i+j-1}$
is equal to
$\lambda_j(i+j+2\alpha)\delta_{i+j+2\alpha,0},$ that is $0.$
So, $\overline{Os}^+(\alpha,0)$ is not only subalgebra, but also an ideal of
$Os^+(\alpha,0).$

Let now $p>0.$ The multiplication on $Os^+(\alpha,\beta,m)$ can be
given by $$x^{(i)}\star x^{(j)}=$$ $${i+j\choose j}x^{(i+j-1)}
+2\beta\delta_{i,0}\delta_{j,0}x^{(p^m-2)}+
2(\alpha\delta_{i,0}\delta_{j,0}
-\beta(\delta_{i,0}\delta_{j,1}+\delta_{i,1}\delta_{j,0}))
x^{(p^m-1)}.$$
In particular,
$$1\star 1=2\beta x^{(p^m-2)}+2\alpha
x^{(p^m-1)},$$
$$1\star x^{(1)}=1-2\beta x^{(p^m-1)},$$
$$1\star x^{(j)}=x^{(j-1)},\quad j>1.$$
If $\alpha=0,$ the algebra
$Os^+(\alpha,\beta,m)$ has the ideal
$$\overline{Os}^+(0,\beta,m)=\{1-2\beta x^{(p^m-1)}, x^{(i)}:
0<i<p^m-1\}$$ of dimension $p^m-1.$

\begin{th} \label{simple} Let
\begin{itemize}
\item $A=Os^+,$ if $p=0,$
\item $A=Os^+(\alpha,\beta),$ if $p=0$ and $\beta\ne 0,$ or
$\beta=0,\alpha\not\in {\frac{1}{2}\bf Z},$
\item $A=\overline{Os}^+(\alpha,0),$ if $p=0$ and $\beta=0,$
$\alpha\in {\frac{1}{2}\bf Z},$
\item $A=Os^+(\alpha,\beta,m),$ if $p>2$ and $\alpha\ne 0,$
\item $A=\overline{Os}^+(0,\beta,m),$ if $p>2$ and
$\alpha=0.$
\end{itemize}
Then $A$ is simple Novikov-Jordan algebra.
\end{th}

{\bf Proof.} Let $U=K[x], p=0.$ Then
$1\star x^i=ix^{i-1}.$ If $J$ is ideal of $Os^+$ and $0\ne X=\sum_{i_0\le i\le i_1}\lambda_ix^i\in J, \lambda_{i_1}\ne 0,$ then
by multiplying $1$ to $X$ $i_1$  times, we obtain the
element $\lambda_{i_1} i_1! x^0\in J.$ Therefore, $1\in J$ and
$x^i=(i+1)^{-1} x^{i+1}\star 1 \in J,$ for any $i\ge 0.$ This means that
$J=Os^+.$

Below in case of $p=0$ we assume that $U$ consists of Laurent
formal power series. Let $X\in A.$ Present $X$ as a linear
combination of basic elements: $X=\sum_{j\ge i}\lambda_jx^j,
\lambda_j\in K,$ if $p=0$ or $X=\sum_{j\ge i}\lambda_jx^{(j)},
\lambda_j\in K,$ if $p>0.$ Set $|X|=i,$ if $\lambda_{i}\ne 0.$

Let $J$ be non-zero ideal of $A.$ Take some $0\ne X\in J.$ Suppose
that $|X|=i.$

{\bf Step 1.} Let $M$ be a set of elements $X\in J,$ such that $|X|\ge 0.$
Prove that $M\ne \emptyset.$

If $p>2$ this statement is trivial.

Let $p=0.$ Suppose that $|X|=i<0.$ There exists some $j>|i|+1>0$ such that
\begin{itemize}
\item $x^j\in A$
\item $j\ne -i-2\alpha,$ and
\item $X'=X\star
x^j=\lambda_i(i+j+2\alpha)x^{i+j-1}+2\beta\lambda_i
x^{i+j-2}+X''\in J,$ where $|X''|\ge i+j-1,$ and $|X''|>i+j-1,$ if
$\beta=0.$
\end{itemize}
If $\beta=0,$ then $\lambda_i(i+j+2\alpha)\ne 0$ and $|X''|>i+j-1.$
Therefore, in this case $X'\ne 0.$ If $\beta\ne 0,$ then
$x^{i+j-2}$ enters with non-zero coefficient $2\beta$ in decomposition of
$X'$ by basic elements. In particular, $X'\ne 0.$
So, in both cases, $X'\ne 0, X'\in J$ and $|X'|\ge 0.$
Therefore, the set
$M=\{X\in J: |X|\ge 0\}$ is not empty.

{\bf Step 2.} Take $0\ne X_0\in M$ with minimal degree $|X_0|.$
Let $|X_0|=i_0\ge 0.$ Prove that $i_0=0,$ if
$(\alpha,\beta)\ne (-1/2,0)$ and $i_0=1,$ if $(\alpha,\beta)=(-1/2,0).$

Consider firstly the case $\alpha=-1/2, \beta=0.$ Suppose that
$i_0>1.$ Then
$1\not\in \overline{Os}^+(-1/2,0),$ but $x^{-1}\in \overline{Os}^+(-1/2,0).$
We have
$$X_1=x^{-1}\star X_0= \lambda_{i_0}(i_0-2)x^{i_0-2}+ X_1'\in J,$$
where $|X_1'|\ge i_0-1.$ Thus, $i_0=2.$ Then
$$X_2=x^2\star X_0=3 \lambda_2  x^3+X_2'\in J, |X_2'|>3.$$
We have
$$X_3=x^{-1}\star X_2=3\lambda_2x+X_3'\in J, |X_3'|>1.$$
Since $\lambda_{i_0}\ne 0, i_0=2,$ we obtain contradiction. So, $i_0=1,$
if $(\alpha,\beta)=(-1/2,0).$

No consider the case, when $(\alpha,\beta)\ne (-1/2,0).$ In this
case $1\in A$ or $1-2\beta x^{(p^m-1)}\in A,$ if $\alpha=0, p>0.$
Suppose that $i_0>0.$ Then $\bar X_0:=1\star X_0\in J$ or $\bar
X_0:=(1-2\beta x^{(p^m-1)})\star X_0\in J$ in case of $\alpha=0,
p>0,$ satisfies the condition
$$0\le |\bar X_0|<|X_0|.$$
Therefore, $\bar X_0=0,$ and $i_0=0.$

So, our statement in step 2  is proved.

{\bf Step 3.} Now prove that $J=A.$

Let $p=0.$ If $\beta\ne 0,$ then $x^{i+2}\in A=Os^+(\alpha,\beta),$
for any $i\in {\bf Z}.$ Thus
$$X_0\star x^{i+2}=2\beta x^{i}+ X_0'\in J, \quad |X_0'|>i,$$
for any $i\in{\bf Z}.$ This means that in $J$ one can get basis
with elements of the form $2\beta x^{i}+X_0', |X_0'|>i, i\in {\bf Z}.$
Recall that in case of $p=0$ our algebras consist of Laurent formal
power series, i.e., infinite sum in positive part is allowed.
Therefore, $J=Os^+(\alpha,\beta),$ if $\beta \ne 0.$

If $p=0,\alpha\ne -1/2, \beta=0,$ then
$$X_0\star x^{i+1}=(i+1+2\alpha)x^i+X_0'\in J,\quad |X_0'|>i.$$
So, $J=Os(\alpha,0),$ if $\alpha\not\in {\frac{1}{2}\bf Z}.$ If
$2\alpha+1\in {\bf Z},$ then $X_0\in\overline{Os}(\alpha,0)$ and
$J=\{x^i: i\in {\bf Z}, i\ne -2\alpha-1\}=\overline{Os}(\alpha,0).$

If $p=0,\alpha= -1/2, \beta=0,$ then
$$X_0\star x^{i}=ix^i+X_0'\in J,\quad |X_0'|>i,$$
for any $i\ne 0.$ So, $J=\overline{Os}^+(-1/2,0).$

If $p>2$ and $x^{(i)}\in A, i>0,$ then $$X_0\star
x^{(i)}=x^{(i-1)}+X_0'\in J, \quad |X_0'|\ge i.$$ So, $J$ has elements
of the form $Y_i=x^{(i)}+\lambda_ix^{(p^m-1)}, 0\le i<p^m-1.$ Notice
that $$Y_1\star Y_1=2 x^{(1)}\in J.$$ Therefore,
$\lambda_1x^{(p^m-1)}\in J.$

Suppose that $J$ has element $x^{(p^m-1)}.$ The multiplication of
this element $p^m-i-1$ times by $x^{(0)}:=1$, gives us the element
$x^{(i)}\in J, i>0.$ Finally,
$1\star x^{(1)}=1-2\beta x^{(p^m-1)}\in J.$
Therefore,
$1\in J,$ and $J=Os(\alpha,\beta,m).$ So, if $\lambda_1\ne 0,$
then our theorem is proved.

Consider the case, $\lambda_1=0.$ Then $x^{(1)}\in J$ and $$Y_0\star
x^{(1)}=1-2\beta x^{(p^m-1)}\in J.$$ So, $\lambda_0=-2\beta.$ Further,
$$Y_0\star Y_0=1\star 1+2\lambda_01\star x^{(p^m-1)}=2\alpha
x^{(p^m-1)}+2\lambda_0x^{(p^m-2)}+2\beta x^{(p^m-2)}\in J.$$ Thus,
$2\alpha x^{(p^m-1)}\in J.$ So, if $\alpha\ne 0,$ then
$J=Os(\alpha,\beta,m).$ If $\alpha=0,$ then
$J=\overline{Os}^+(0,\beta,m).$ Theorem is proved completely.

{\bf Remark.} If in the case $p=0$ we consider Laurent polynomials
algebra instead of formal power series, then theorem \ref{simple}
is not true. Below we give one counterexample.

Set $\overline{x^i}=x^i,$ if  $ i<-1$
and $\overline{x^i}= x^i+2(i+1)^{-1}\beta x^{i-1},$ if $i>-1.$
Let  $\overline{Os}(0,\beta)$  be a subspace of $Os(0,\beta)$ generated by
elements $\overline{x^i}, i\ne -1.$
Notice that, if $\alpha=0, \beta=0$  these two definitions of subspace
$\overline{Os}(0,0)$ are coincide.

\begin{lm} \label{22nov}$(p=0).$ Let $U=K[x^{\pm 1}].$
The subspace $\overline{Os}(0,\beta)\subset Os(0,\beta)$ is close
under multiplication  $(a,b)\mapsto a\star b=\der(a b)+2\beta
x^{-2} ab.$ Moreover, it is an ideal of codimension $1.$
\end{lm}

{\bf Proof.} It is obvious.

 {\bf Remark.} The algebra $\overline{Os}(0,0)=\{x^i:
i\ne -1\}$ has a nontrivial central extension. The cocycle $\psi$
corresponding to this central extension can be given by
$$\psi(x^i,x^j)=\delta_{i+j,0}.$$
It has the following property:
$$\psi(\{x^i,x^j\},x^s)+\psi(\{x^j,x^s\},x^i)+\psi(\{x^s,x^i\},x^j)=2,$$
if
$$\psi(\{x^i,x^j\},x^s)+\psi(\{x^j,x^s\},x^i)+\psi(\{x^s,x^i\},x^j)
\ne 0,$$
i.e., if $i+j+s=1.$ In the case of  characteristic $p>0$ the definition of
$\psi$ should be slightly changed:
$$\psi(x^{(i)},x^{(j)})=\frac{1}{p}{p^m\choose
i}\delta_{i+j,p^m}.$$ Notice that ${p^m\choose i}\equiv
0(mod\,p),$ if $0\le i<p^m,$ and this definition is correct.

\section{Leibniz dual algebras}

An algebra $A$ with multiplication $\circ$ is called {\it (left)
Leibniz,} if it satisfies the identity
$$(a\circ b)\circ c-a\circ (b\circ c)+b\circ (a\circ c)=0.$$
Similarly, an algebra with identity
$$a\circ (b\circ c)-(a\circ b)\circ c+(a\circ c)\circ b=0$$
is called right Leibniz \cite{loday}. An algebra $(A,\times)$ is
called (left) Leibniz dual, if it satisfies the identity
$$(a\times b)\times c=a\times (b\times c)+a\times (c\times b).$$

\begin{prp} Let $A$ be a Leibniz dual algebra. Then it satisfies
(right) Tortken identity.
\end{prp}

{\bf Proof.} According (\ref{Leibnizdual})
 $$(a,b,c)=a\times(b\times c)-(a\times b)\times c= -a\times(c\times b).$$
According (\ref{rcommut}) $$(a,b,c)\times d=-(a\times d)\times
(c\times b).$$ By analogous reasons $$-(a,d,c)\times b=(a\times
b)\times (c\times b).$$ Thus $$(a\times b)\times (c\times
d)-(a\times d)\times(c\times b)=(a,b,c)\times d-(a,d,c)\times c.$$

{\bf Example 1.} In \cite{Lodaycupprod} it is proved that the
Leibniz cohomology group of trivial module $H_{lei}^*(L,K)$ for
Leibniz algebra is Leibniz dual under the cup product. Statements
of this paper can be generalized. Let us give some of such
generalizations.

\begin{prp} \label{loday1}
Let $\frak g$ be a right Leibniz algebra and $R$ be a left
Leibniz dual algebra. Then the tensor product ${\frak g}\otimes R$
equipped with the multiplication
$$(x\otimes r)\circ (y\otimes s) =[x,y]\otimes rs$$
is a right-symmetric algebra.
\end{prp}

{\bf Proof.} By definition
$$(x\otimes r)\circ ((y\otimes s)\circ (z\otimes t))-
((x\otimes r)\circ (y\otimes s))\circ (z\otimes t))$$
$$-(x\otimes r)\circ ((z\otimes t)\circ (y\otimes s))+
((x\otimes r)\circ (z\otimes t))\circ (y\otimes s))=$$

$$[x,[y,z]]\otimes r(st)-\mathop{[[x,y],z]\otimes (rs)t}\limline
-[x,[z,y]]\otimes r(ts)+\mathop{[[x,z],y]\otimes (rt)s}\limsim=$$
(by left Leibniz dual condition)
$$\mathop{[x,[y,z]]\otimes r(st)}\limline
-\mathop{[[x,y],z]\otimes r(st)}\limline-
\mathop{[[x,y],z]\otimes r(ts)}\limsim$$
$$-\mathop{[x,[z,y]]\otimes r(ts)}\limsim
+\mathop{[[x,z],y]\otimes r(ts)}\limsim+
\mathop{[[x,z],y]\otimes r(st)}\limline=$$
(by right  Leibniz dual condition)
$$=0.$$

\begin{crl} {\rm (Proposition $1.3$ of \cite{Lodaycupprod})}
 Let $\frak g$ be a right Leibniz algebra and $R$ be a left
Leibniz dual algebra. Then the tensor product ${\frak g}\otimes R$
equipped with the multiplication
$$(x\otimes r)\circ (y\otimes s) =[x,y]\otimes rs-[y,x]\otimes sr$$
is a Lie algebra.
\end{crl}

Second generalization concerns cup-products.
Cup-products of Leibniz cohomologies  can be defined not
only for pairing of trivial modules as it done in
\cite{Lodaycupprod}. Let $\frak g$ be a left  Leibniz algebra, $M$ be a
symmetric module, $N$ be an anti-symmetric module and $S$ be an
anti-symmetric module. Suppose that there is given a bilinear map
$$M\times N\rightarrow S,\quad (m,n)\mapsto m\cup n$$ such that $$
X(m\cup n)=X(m)\cup n+m\cup X(n),$$ for any $X\in {\frak g}, m\in
M, n\in N.$ This cup product can be prolonged until the bilinear
map
$$C^k_{lei}({\frak g},M)\times C_{lei}^l({\frak g},N)\rightarrow
C^{k+l}_{lei}({\frak g},S)$$ by
$$\psi\cup\phi(X_1,\ldots,X_{k+l})=$$ $$\sum_{\sigma
}sign\,\sigma\,\psi(X_{\sigma(1)},\ldots,X_{\sigma(k)})\cup
\phi(X_{\sigma(k+1)},\ldots,X_{\sigma(k+l)}).$$ Summation here is
under $\sigma\in Sym_{k+l},$ such that $\sigma(1)<\cdots\sigma(k)$
and  $\sigma(k+1)<\cdots<\sigma(k+l), \sigma(k+l)=k+l.$ This
bilinear map is compatible with action of a coboundary operator
$$d(\psi\cup\phi)=d\psi\cup\phi+(-1)^{|\psi|}\psi\cup d\phi.$$ If
$M=K$ be the trivial module and $N$ is an antisymmetric module
and $(K,N) \rightarrow N, (\lambda,n)\mapsto \lambda n$ be a natural
bilinear map, then the the cup-product satisfies the
Leibniz dual rule $$\psi\cup (\phi\cup
\chi)=(\psi\cup\phi+(-1)^{|\psi||\phi|}\phi\cup\psi)\cup\chi,$$
where $\psi\in C^*_{lei}({\frak g},K), \phi\in C^*_{lei}({\frak g},N), \chi\in
C^*_{lei}({\frak g},N).$ This means that $H^*_{lei}({\frak g},K)$ is not only
right Leibniz dual algebra, but also $H^*_{lei}({\frak g},N)$ is a
left module over right Leibniz dual algebra $H^*_{lei}({\frak g},K),$
if $N$ is antisymmetric.

Another generalization concerns cup-products of right-symmetric algebras.

{\bf Example 2.} Let $A$ be a right-symmetric algebra and $M$ be a
right-symmetric module with a cup product $$M\times M\rightarrow
M, (m,n)\mapsto m\cup n.$$ For instance we can take as $M$  the
trivial module  with usual multiplication $K\times K\rightarrow K,
(\lambda,\mu)\mapsto \lambda\mu.$ Let $$C^*_{rsym}(A,M)\times
C^*_{rsym}(A,M)\rightarrow C^*_{rsym}(A,M), \quad
(\psi,\phi)\mapsto \psi\cup\phi$$ be the cup product of
corresponding cochain complexes \cite{DzhumaKos}. Then this cup
product satisfies Leibniz dual law $$(\psi\cup \phi)\cup
\chi=\psi\cup(\phi\cup \chi+(-1)^{|\phi||\chi|}\chi\cup\phi),$$
where $\psi\in C^*_{rsym}(A,M), \phi\in C^*_{rsym}(A,M), \chi\in
C^*_{rsym}(A,M).$ For instance, the right-symmetric cohomology of
the trivial module $H_{rsym}^*(A,K)$ is a left Leibniz dual algebra under
the cup product and $H^*_{rsym}(A,M)$ is a right module over left Leibniz
dual algebra $H^*_{rsym}(A,K).$

Proof of this statement repeats arguments of \cite{Lodaycupprod}.

{\bf Example 3.} Let $A =K[x]$ and
$a\star b=a\int_0^x b\,dx.$ Then $(A,\star)$ is a left Leibniz dual.
It is easy to see
that identity (\ref{Leibnizdual}) is equivalent to the
condition of integration by parts.

\section{Tortken algebras in characteristic $p>0$}

\begin{th}\label{tirtforcharp}
Let $K$ be a field of characteristic $char\, K=p>0,$ $(A,\circ)$
be an associative commutative algebra and $D$ be a derivation of
$(A,\circ).$ Let $\square=\square_{k,l} $  be a new multiplication
on $A$ depending on some integers $0\le k\le l,$ such that
$$a\square b=D(D^{p^k-1}(a)\circ D^{p^l-1}(b)
+D^{p^l-1}(a)\circ D^{p^k-1}(b)),$$
if $k\ne l$ and
$$a\square b=D(D^{p^k-1}(a)\circ D^{p^k-1}(b)),$$
if $k=l.$
Then

{\rm i)} $(A,\square)$ is Tortken, if $k=l, p>0$ or $l=k+1, p=2.$

{\rm ii)} If $k\ne l, p>2$ or $l-k>1, p=2,$ then there exists some associative
commutative algebra $(A,\circ),$ such that $(A,\square)$ is not Tortken.
\end{th}

{\bf Proof.} Define multiplications $\star$ and $\star'$ by
$$a\star b= D(D^{p^k-1}(a)\circ D^{p^l-1}(b)),$$
$$a\star' b= D(D^{p^l-1}(a)\circ D^{p^k-1}(b)).$$
Then
$$a\square b= a\star b+a\star' b$$
if $k\ne l$ and
$$a\square b=a\star b,$$
if $k=l.$

Let
$$f(a,b,c,d)=$$
$$(a\square b)\square(c\square d)-(a\square d)\square(c\square b)$$
$$-(a,b,c)^{\square}\square d+(a,d,c)^{\square}\square b.$$
For the multiplication $\star$ we have
$$(a,b,c)^{\star}=a\star (b\star c)-(a\star b)\star c=$$

$$D(D^{p^k-1}(a)\circ D^{p^l-1}(D (D^{p^k-1}(b)\circ D^{p^l-1}(c))))$$
$$-D(D^{p^k-1}(D(D^{p^k-1}(a)\circ D^{p^l-1}(b)))\circ D^{p^l-1}(c))=$$

$$D(D^{p^k-1}(a)\circ D^{p^l}(D^{p^k-1}(b)\circ D^{p^l-1}(c))-
D^{p^k}(D^{p^k-1}(a)\circ D^{p^l-1}(b))\circ D^{p^l-1}(c))=$$

$$D\left(\mathop{D^{p^k-1}(a)\circ D^{p^k+p^l-1}(b)\circ D^{p^l-1}(c)}
\limline
+D^{p^k-1}(a)\circ D^{p^k-1}(b)\circ D^{2p^l-1}(c)\right.$$
$$\left.-D^{2p^k-1}(a)\circ D^{p^l-1}(b)\circ D^{p^l-1}(c)
-\mathop{D^{p^k-1}(a)\circ D^{p^k+p^l-1}(b)\circ D^{p^l-1}(c)}\limline
\right)=$$

$$D\left(D^{p^k-1}(a)\circ D^{p^k-1}(b)\circ D^{2p^l-1}(c)
-D^{2p^k-1}(a)\circ D^{p^l-1}(b)\circ D^{p^l-1}(c)\right).
$$
Thus
$$(a,b,c)^{\star}\star d=$$

$$D(D^{p^k}\left(D^{p^k-1}(a)\circ D^{p^k-1}(b)\circ D^{2p^l-1}(c)\right.$$
$$
\left.-D^{2p^k-1}(a)\circ D^{p^l-1}(b)\circ D^{p^l-1}(c)\right) \circ D^{p^l-1}(d))=
$$

$$
D\left(
D^{2p^k-1}(a)\circ D^{p^k-1}(b)\circ D^{2p^l-1}(c)\circ D^{p^l-1}(d)\right.$$
$$+D^{p^k-1}(a)\circ D^{2p^k-1}(b)\circ D^{2p^l-1}(c)\circ D^{p^l-1}(d)$$
$$+D^{p^k-1}(a)\circ D^{p^k-1}(b)\circ D^{p^k+2p^l-1}(c)\circ D^{p^l-1}(d)$$
$$\mathop{\left.-D^{3p^k-1}(a)\circ D^{p^l-1}(b)\circ D^{p^l-1}(c)\right)
\circ D^{p^l-1}(d)}\limline$$
$$\left.-D^{2p^k-1}(a)\circ D^{p^k+p^l-1}(b)\circ
D^{p^l-1}(c)\right) \circ D^{p^l-1}(d)$$
$$\left.\mathop{-D^{2p^k-1}(a)\circ D^{p^l-1}(b)\circ
D^{p^k+p^l-1}(c)\circ D^{p^l-1}(d)}\limeq\right).$$

Similarly,
$$(a,d,c)^{\star}\star b=$$

$$D\left(
D^{2p^k-1}(a)\circ D^{p^k-1}(d)\circ D^{2p^l-1}(c)\circ D^{p^l-1}(b)\right.$$
$$+D^{p^k-1}(a)\circ D^{2p^k-1}(d)\circ D^{2p^l-1}(c)\circ D^{p^l-1}(b)$$
$$+D^{p^k-1}(a)\circ D^{p^k-1}(d)\circ D^{p^k+2p^l-1}(c)\circ D^{p^l-1}(b)$$
$$\mathop{\left.-D^{3p^k-1}(a)\circ D^{p^l-1}(d)\circ D^{p^l-1}(c)\right)
\circ D^{p^l-1}(b)}\limline$$
$$-D^{2p^k-1}(a)\circ D^{p^k+p^l-1}(d)\circ
D^{p^l-1}(c) \circ D^{p^l-1}(d)$$
$$\left.\mathop{-D^{2p^k-1}(a)\circ D^{p^l-1}(d)\circ
D^{p^k+p^l-1}(c) \circ D^{p^l-1}(b)}\limeq\right).$$
Therefore,
$$(a,b,c)^{\star}\star d-(a,d,c)^{\star}\star b=$$

$$
D\left(
D^{2p^k-1}(a)\circ D^{p^k-1}(b)\circ D^{2p^l-1}(c)\circ D^{p^l-1}(d)\right.$$
$$+D^{p^k-1}(a)\circ D^{2p^k-1}(b)\circ D^{2p^l-1}(c)\circ D^{p^l-1}(d)$$
$$+D^{p^k-1}(a)\circ D^{p^k-1}(b)\circ D^{p^k+2p^l-1}(c)\circ D^{p^l-1}(d)$$
$$-D^{2p^k-1}(a)\circ D^{p^k+p^l-1}(b)\circ
D^{p^l-1}(c) \circ D^{p^l-1}(d)$$
$$-D^{2p^k-1}(a)\circ D^{p^k-1}(d)\circ D^{2p^l-1}(c)\circ D^{p^l-1}(b)$$
$$-D^{p^k-1}(a)\circ D^{2p^k-1}(d)\circ D^{2p^l-1}(c)\circ D^{p^l-1}(b)$$
$$-D^{p^k-1}(a)\circ D^{p^k-1}(d)\circ D^{p^k+2p^l-1}(c)\circ D^{p^l-1}(b)$$
$$\left.+D^{2p^k-1}(a)\circ D^{p^k+p^l-1}(d)\circ
D^{p^l-1}(c) \circ D^{p^l-1}(b)\right).$$
On the other hand,
$$(a\star b)\star (c\star d)=$$
$$D\left(D^{p^k}(D^{p^k-1}(a)\circ D^{p^l-1}(b))
\circ D^{p^l}(D^{p^k-1}(c)\circ D^{p^l-1}(d))\right)=$$

$$D\left(\mathop{D^{2p^k-1}(a)\circ D^{p^l-1}(b)
\circ D^{p^k+p^l-1}(c)\circ D^{p^l-1}(d)}\limline\right.$$
$$+
D^{p^k-1}(a)\circ D^{p^k+p^l-1}(b)
\circ D^{p^k+p^l-1}(c)\circ D^{p^l-1}(d)$$
$$
+D^{2p^k-1}(a)\circ D^{p^l-1}(b)
\circ D^{p^k-1}(c)\circ D^{2p^l-1}(d)$$ $$\left.+
D^{p^k-1}(a)\circ D^{p^k+p^l-1}(b)
\circ D^{p^k-1}(c)\circ D^{2p^l-1}(d)\right).$$
Similarly,
$$(a\star d)\star (c\star b)=$$

$$D\left(\mathop{D^{2p^k-1}(a)\circ D^{p^l-1}(d)
\circ D^{p^k+p^l-1}(c)\circ D^{p^l-1}(b)}\limline\right.$$
$$+
D^{p^k-1}(a)\circ D^{p^k+p^l-1}(d)
\circ D^{p^k+p^l-1}(c)\circ D^{p^l-1}(b)$$
$$
+D^{2p^k-1}(a)\circ D^{p^l-1}(d)
\circ D^{p^k-1}(c)\circ D^{2p^l-1}(b)$$ $$\left.+
D^{p^k-1}(a)\circ D^{p^k+p^l-1}(d)
\circ D^{p^k-1}(c)\circ D^{2p^l-1}(b)\right).$$
Thus,
$$(a\star b)\star (c\star d)-(a\star d)\star (c\star b)=$$

$$D\left(
D^{p^k-1}(a)\circ D^{p^k+p^l-1}(b)
\circ D^{p^k+p^l-1}(c)\circ D^{p^l-1}(d)\right.$$
$$
+D^{2p^k-1}(a)\circ D^{p^l-1}(b)
\circ D^{p^k-1}(c)\circ D^{2p^l-1}(d)$$ $$+
D^{p^k-1}(a)\circ D^{p^k+p^l-1}(b)
\circ D^{p^k-1}(c)\circ D^{2p^l-1}(d)$$

$$-D^{p^k-1}(a)\circ D^{p^k+p^l-1}(d)
\circ D^{p^k+p^l-1}(c)\circ D^{p^l-1}(b)$$
$$
-D^{2p^k-1}(a)\circ D^{p^l-1}(d)
\circ D^{p^k-1}(c)\circ D^{2p^l-1}(b)$$ $$\left.-
D^{p^k-1}(a)\circ D^{p^k+p^l-1}(d)
\circ D^{p^k-1}(c)\circ D^{2p^l-1}(b)\right).$$

Now we are ready to establish i). We see that, if $k=l,$ then
$$f(a,b,c,d)=$$

$$D\left(
\mathop{D^{p^k-1}(a)\circ D^{2p^k-1}(b)\circ D^{2p^k-1}(c)\circ D^{p^k-1}(d)}\limline\right.$$
$$-\mathop{D^{2p^k-1}(a)\circ D^{2p^k-1}(b)\circ
D^{p^k-1}(c) \circ D^{p^k-1}(d)}\limeq$$
$$-\mathop{D^{p^k-1}(a)\circ D^{2p^k-1}(d)\circ D^{2p^k-1}(c)\circ
D^{p^k-1}(b)}\limsim$$
$$+\mathop{D^{2p^k-1}(a)\circ D^{2p^k-1}(d)\circ
D^{p^k-1}(c) \circ D^{p^k-1}(b)}\limsimeq$$

$$-\mathop{D^{p^k-1}(a)\circ D^{2p^k-1}(b)
\circ D^{2p^k-1}(c)\circ D^{p^k-1}(d)}\limline$$
$$-\mathop{D^{2p^k-1}(a)\circ D^{p^k-1}(b)
\circ D^{p^k-1}(c)\circ D^{2p^k-1}(d)}\limsimeq$$
$$+\mathop{D^{p^k-1}(a)\circ D^{2p^k-1}(d)
\circ D^{2p^k-1}(c)\circ D^{p^k-1}(b)}\limsim$$
$$\left.+\mathop{D^{2p^k-1}(a)\circ D^{p^k-1}(d)
\circ D^{p^k-1}(c)\circ D^{2p^k-1}(b)}\limeq
\right)$$

$$=0.$$

The case $p=2, l=k+1$ can be established by analogous way. One
just need to use the following relations $2p^k=p^l,
p^k+p^l-1=3p^k-1.$

To prove ii) take the algebra of divided powers $O_1(m)$ as $A.$
Then, substitutions $a=x^{(p^k-1)}, b= x^{(2p^k-1)}, c=
x^{(p^l-1)}, d=x^{(2p^l)}$ in formulas for $f(a,b,c,d)$ show that
$$f(x^{(p^k-1)}, x^{(2p^k-1)}, x^{(p^l-1)}, x^{(2p^l)})=-1,$$
if $l>k, p>2$ or $l-k>1, p=2.$

{\bf Remark.} We was not able to construct in $O_1(m)$ some Novikov
multiplication $\ast,$ such that $a\square b= a\ast b+b\ast a.$
We do not know whether the divided power algebra $O_1(m)$
under Tortken multiplication $\square$ constructed in
theorem~\ref{tirtforcharp} is special.

\section{Consequences of Tortken identity of degree 5}

\begin{th}\label{oneki} Let $A$ be a commutative Tortken algebra
and $a,b,c,x,y$ are any elements of $ A.$ Then

{\rm i)} $$(a,y,b)\star(x\star c)+(b,y,c)\star (x\star a)+(c,y,a)\star (x\star b)$$
$$-(a\star y)\star (b,x,c)-(b\star y)\star (c,x,a)-(c\star y)\star (a,x,b)=0,$$

{\rm ii)}  $(p\ne 3)$
$$((x\star a)\star (y\star b)-(x\star b)\star (y\star a))\star c$$
$$+((x\star b)\star (y\star c)-(x\star c)\star (y\star b))\star a$$
$$+((x\star c)\star (y\star a)-(x\star a)\star (y\star c))\star b=0,$$

{\rm iii)}
$$(x,y\star a,b)\star c+(x,y\star b,c)\star a+(x,y\star c,a)\star b$$
$$-(x,y\star b,a)\star c-(x,y\star c,b)\star a-(x,y\star a,c)\star b=0,$$

{\rm iv)}
$$
(((x\star a)\star b)\star c)\star y
+(((x\star b)\star c)\star a)\star y
+(((x\star c)\star a)\star b)\star y
$$
$$-(((x\star b)\star a)\star c)\star y
-(((x\star c)\star b)\star a)\star y
-(((x\star a)\star c)\star b)\star y
=0.$$
\end{th}

{\bf Proof.}
i) According (\ref{main}) for $a_1=a,a_2=b,a_3=c,$
$$(a_{\sigma(1)},x\star a_{\sigma(3)},a_{\sigma(2)})\star y=
$$
$$(a_{\sigma(1)},y,a_{\sigma(2)})\star (x\star a_{\sigma(3)})$$
$$-(a_{\sigma(1)}\star y)\star (a_{\sigma(2)}\star (x\star
a_{\sigma(3)}))+(a_{\sigma(1)}\star(x\star a_{\sigma(3)}))\star
(a_{\sigma(2)}\star y).$$ Let $Sym_3^0$ be the cyclic subgroup of
$Sym_3$ of order 3. Therefore,
$$\sum_{\sigma\in
Sym_3^0}(a_{\sigma(1)},y,a_{\sigma(2)})\star(x\star
a_{\sigma(3)})$$ $$-(a_{\sigma(1)}\star y)\star
(a_{\sigma(2)}\star (x\star
a_{\sigma(3)}))+(a_{\sigma(1)}\star(x\star a_{\sigma(3)}))\star
(a_{\sigma(2)}\star y) $$ $$=y\star \sum_{\sigma\in
Sym_3^0}(a_{\sigma(1)},x\star a_{\sigma(3)},a_{\sigma(2)})=$$ (by
corollary~\ref{on}) $$=0.$$ Since
$$(a,y,b)\star(x\star c)+(b,y,c)\star (x\star a)+(c,y,a)\star (x\star b)$$
$$-\mathop{(a\star y)\star (b\star (x\star c))}\limline+
\mathop{(b\star y)\star (a\star (x\star c))}\limeq$$
$$-\mathop{(b\star y)\star (c\star (x\star a))}\limeq+
\mathop{(c\star y)\star (b\star (x\star a))}\limsim$$
$$-\mathop{(c\star y)\star (a\star (x\star b))}\limsim+\mathop{
(a\star y)\star (c\star (x\star b))}\limline=$$

$$(a,y,b)\star(x\star c)+(b,y,c)\star (x\star a)+(c,y,a)\star (x\star b)$$
$$-(a\star y)\star (b,x,c)-(b\star y)\star (c,x,a)-(c\star y)\star (a,x,b),$$
our statement is proved.

ii) Set
$$g=((x\star a)\star (y\star b)-(x\star b)\star (y\star a))\star c$$
$$+((x\star b)\star (y\star c)-(x\star c)\star (y\star b))\star a$$
$$+((x\star c)\star (y\star a)-(x\star a)\star (y\star c))\star b.$$
According (\ref{main}),
$$((x,a,y)\star b-(x, b,y)\star a)\star c$$
$$+((x,b,y)\star c-(x,c,y)\star b)\star a$$
$$+((x,c,y)\star a-(x,a,y)\star c)\star b=$$

$$((x,a,y)\star b)\star c-((x,a,y)\star c)\star b$$
$$+((x,b,y)\star c)\star a-((x,b,y)\star a)\star c$$
$$-((x,c,y)\star b)\star a+((x,c,y)\star a)\star b=$$

$$(a,(x,b,y),c)+(b,(x,c,y),a)+(c,(x,a,y),b)=$$

$$(a,x\star (y\star b),c)-(a,y\star (x\star b),c)$$
$$(b,x\star (y\star c),a)-(b,y\star (x\star c),a)$$
$$(c,x\star (y\star a),b)-(c,y\star (x\star a),b).$$
We have
$$(a,x\star (y\star b),c)=$$
(by corollary \ref{on})
$$(a,c\star x,y\star b)-(c,a\star x,y\star b).$$
Similarly,
$$(b,x\star (y\star c),a)=
(b,a\star x,y\star c)-(a,b\star x,y\star c)
$$
$$(c,x\star (y\star a),b)=
(c,b\star x,y\star a)-(b,c\star x,y\star a)
$$

$$-(a,y\star (x\star b),c)=
-(a,c\star y,x\star b)+(c,a\star y,x\star b),$$

$$-(b,y\star (x\star c),a)=
-(b,a\star y,x\star c)+(a,b\star y,x\star c),$$

$$-(c,y\star (x\star a),b)=
-(c,b\star y,x\star a)+(b,c\star y,x\star a).$$
Hence
$$g=$$
$$(a,c\star x,y\star b)-(c,a\star x,y\star b)$$
$$+(b,a\star x,y\star c)-(a,b\star x,y\star c)$$
$$+(c,b\star x,y\star a)-(b,c\star x,y\star a)$$
$$-(a,c\star y,x\star b)+(c,a\star y,x\star b)$$
$$-(b,a\star y,x\star c)+(a,b\star y,x\star c)$$
$$-(c,b\star y,x\star a)+(b,c\star y,x\star a).$$
Collect together all associators with the first element $a:$
$$(a,c\star x,y\star b)-(a,b\star x,y\star c)
-(a,c\star y,x\star b)+(a,b\star y,x\star c)=$$

$$a\star((c\star x)\star(y\star b))-(a\star (c\star x))\star(y\star b)$$
$$-a\star((b\star x)\star(y\star c))+(a\star(b\star x))\star (y\star c)$$
$$-a\star((c\star y)\star(x\star b))+(a\star (c\star y))\star (x\star b)$$
$$+a\star((b\star y)\star(x\star c))-(a\star(b\star y))\star(x\star c)=$$

$$2a\star((c\star x)\star(y\star b))-(a\star (c\star x))\star(y\star b)$$
$$-2a\star((b\star x)\star(y\star c))+(a\star(b\star x))\star (y\star c)$$
$$+(a\star (c\star y))\star (x\star b)-(a\star(b\star y))\star(x\star c).$$
Similarly,
$$(b,a\star x,y\star c)-(b,c\star x,y\star a)
-(b,a\star y,x\star c)+(b,c\star y,x\star a)=$$

$$2b\star((a\star x)\star(y\star c))-(b\star (a\star x))\star(y\star c)$$
$$-2b\star((c\star x)\star(y\star a))+(b\star(c\star x))\star (y\star a)$$
$$+(b\star (a\star y))\star (x\star c)-(b\star(c\star y))\star(x\star a)$$
and
$$(c,b\star x,y\star a)-(c,b\star x,y\star b)
-(c,b\star y,x\star a)+(c,a\star y,x\star b)=$$

$$2c\star((b\star x)\star(y\star a))-(c\star (b\star x))\star(y\star a)$$
$$-2c\star((a\star x)\star(y\star b))+(c\star(a\star x))\star (y\star b)$$
$$+(c\star (b\star y))\star (x\star a)-(c\star(a\star y))\star(x\star b).$$
Hence
$$g=$$

$$2a\star((c\star x)\star(y\star b))-(a\star (c\star x))\star(y\star b)$$
$$-2a\star((b\star x)\star(y\star c))+(a\star(b\star x))\star (y\star c)$$
$$+(a\star (c\star y))\star (x\star b)-(a\star(b\star y))\star(x\star c)$$
$$+2b\star((a\star x)\star(y\star c))-(b\star (a\star x))\star(y\star c)$$
$$-2b\star((c\star x)\star(y\star a))+(b\star(c\star x))\star (y\star a)$$
$$+(b\star (a\star y))\star (x\star c)-(b\star(c\star y))\star(x\star a)$$
$$+2c\star((b\star x)\star(y\star a))-(c\star (b\star x))\star(y\star a)$$
$$-2c\star((a\star x)\star(y\star b))+(c\star(a\star x))\star (y\star b)$$
$$+(c\star (b\star y))\star (x\star a)-(c\star(a\star y))\star(x\star b)=$$

$$-2g$$
$$-(a\star (c\star x))\star(y\star b)+(a\star(b\star x))\star (y\star c)$$
$$+(a\star (c\star y))\star (x\star b)-(a\star(b\star y))\star(x\star c)$$
$$-(b\star (a\star x))\star(y\star c)+(b\star(c\star x))\star (y\star a)$$
$$+(b\star (a\star y))\star (x\star c)-(b\star(c\star y))\star(x\star a)$$
$$-(c\star (b\star x))\star(y\star a)+(c\star(a\star x))\star (y\star b)$$
$$+(c\star (b\star y))\star (x\star a)-(c\star(a\star y))\star(x\star b).$$
Collect in the last expression all terms of the form $(\cdots)\star (y\star b).$
We obtain
$$-(a\star (c\star x))\star(y\star b)+(c\star(a\star x))\star (y\star b)=$$
$$(c,x,a)\star (y\star b).$$
Collect similarly, all terms of the form $(\cdots)\star(u\star v),$ where
$u=x,y,$ $v=a,b,c.$ We have
$$g=$$
$$-2g+(a,x,b)\star (y\star c)+(b,x,c)\star (y\star a)+(c,x,a)\star (y\star b)$$
$$-(a,y,b)\star (x\star c)-(b,y,c)\star (x\star a)-(c,y,a)\star (x\star b)=$$
(by (i))
$$-2 g.$$
Thus, $3g=0,$ and $g=0,$ if $p\ne 3.$

iii) Follows from corollaries~\ref{on} and \ref{onon} .

iv) Follows from proposition~\ref{rsym3}.

\begin{crl}\label{7november} If $A$ is commutative algebra with identity
$(\ref{main})$ over a field of characteristic $p\ne 3,$ then
$$(a,(b,x,c),y)+(b,(c,x,a),y)+(c,(a,x,b),y)$$
$$-(a,(b,y,c),x)-(b,(c,y,a),x)-(c,(a,y,b),x)=0,$$
for any $a,b,c,x,y\in A.$
\end{crl}

{\bf Proof.} By commutativity rule,
$$((x\star a)\star (y\star b)-(x\star b)\star (y\star a))\star c$$
$$+((x\star b)\star (y\star c)-(x\star c)\star (y\star b))\star a$$
$$+((x\star c)\star (y\star a)-(x\star a)\star (y\star c))\star b=$$

$$((a\star x)\star (b\star y)-(a\star y)\star(b\star x))\star c$$
$$+((b\star x)\star (c\star y)- (b\star y)\star (c\star x))\star a$$
$$+((c\star x)\star (a\star y)- (c\star y)\star (a\star x))\star b=$$
(identity (\ref{main}))
$$((a,x,b)\star y-(a,y,b)\star x)\star c$$
$$+((b,x,c)\star y- (b,y,c)\star x)\star a$$
$$+((c,x,a)\star y- (c,y,a)\star x)\star b=$$

$$((a,x,b)\star y)\star c-((a,y,b)\star x)\star c$$
$$+((b,x,c)\star y)\star a-((b,y,c)\star x)\star a$$
$$+((c,x,a)\star y)\star b-((c,y,a)\star x)\star b=$$
(by corollary \ref{onon})

$$((a,x,b)\star y)\star c-((a,y,b)\star x)\star c$$
$$+((b,x,c)\star y)\star a-((b,y,c)\star x)\star a$$
$$+((c,x,a)\star y)\star b-((c,y,a)\star x)\star b$$

$$-((a,x,b)\star c)\star y+((a,y,b)\star c)\star x$$
$$-((b,x,c)\star a)\star y+((b,y,c)\star a)\star x$$
$$-((c,x,a)\star b)\star y+((c,y,a)\star b)\star x=$$

$$+(a,(b,x,c),y)+(b,(c,x,a),y)+(c,(a,x,b),y)$$
$$-(a,(b,y,c),x)-(b,(c,y,a),x)-(c,(a,y,b),x).$$
It remains to use theorem~\ref{oneki}.

{\bf Remark.} Very likely that all identities of degree 5 of
free special Tortken algebra follow from four types of
identities mentioned in theorem \ref{oneki}, if $p>3.$
It will be interesting to find some $s$-identity (if such identity
exists), i.e., an identity that is true for any special Tortken
algebra, but is not true for some commutative Tortken algebra.

\bigskip
\bigskip
\begin{center}
{\em ACKNOWLEDGEMENTS}
\end{center}
\nopagebreak
I am thankful to referee for careful reading the text and
essential remarks. I am grateful to INTAS foundation and
Royal Swedish Academy of Sciences for support.

\end{document}